%% file: Estimates.tex
\documentclass[preprint,3p]{elsarticle}
\pdfoutput=1

\usepackage{amsmath,amssymb}
\usepackage{subfig}
\usepackage{mathtools}
\usepackage{algorithm}
\usepackage{algpseudocode}

\newcommand{\ds}{\displaystyle}

\journal{Journal of Computational and Applied Mathematics}
\biboptions{comma,sort&compress}

\begin{document}
\begin{frontmatter}
\title{Computational complexity and memory usage for multi-frontal direct solvers in structured mesh finite elements}
\author[kaust]{Nathan~Collier}
\author[bilbao]{David~Pardo\corref{cor}}
\ead{}
\author[agh]{Maciej~Paszynski}
\author[kaust]{Victor~M.~Calo}
\address[kaust]{Applied Mathematics and Computational Science and Earth and Environmental Sciences and Engineering\\
King Abdullah University of Science and Technology\\Thuwal, Saudi Arabia}
\address[bilbao]{Department of Applied Mathematics, Statistics, and
  Operational Research\\The University of the Basque Country and Ikerbasque\\Bilbao, Spain}
\address[agh]{AGH University of Science and Technology\\Department of Computer Science\\Krakow, Poland}
\cortext[cor]{Corresponding author}
\begin{abstract}
  The multi-frontal direct solver is the state-of-the-art algorithm
  for the direct solution of sparse linear systems. This paper
  provides computational complexity and memory usage estimates for the
  application of the multi-frontal direct solver algorithm on linear
  systems resulting from B-spline-based isogeometric finite elements,
  where the mesh is a structured grid. Specifically we provide the
  estimates for systems resulting from $C^{p-1}$ polynomial B-spline
  spaces and compare them to those obtained using $C^0$ spaces.
\end{abstract}
\end{frontmatter}
\section{Introduction}\input{intro.tex}

\section{Multi-frontal direct solver algorithm}\input{alg.tex}

\section{Computational complexity and memory usage}\input{complexity.tex}
\section{Numerical Results}\input{results.tex}
\section{Conclusions}\input{conc.tex}
\section{Acknowledgements}
DP has been partially supported by the Spanish Ministry of Sciences
and Innovation Grant MTM2010-16511. MRP has been partially supported
by the Polish MNiSW grant no. NN 519 405737 and NN519 447 739.
\bibliographystyle{elsarticle-num}
\bibliography{paper.bib}
\end{document}

%% file: intro.tex
The main purpose of this paper is to explain in detail how one can
derive estimates for computational complexity and memory usage of the
multifrontal direct solver algorithm as applied to B-spline-based
\cite{Farin2002,Piegl1995} isogeometric finite elements
\cite{Hughes2005,Cottrell2006}. In this paper we generalize results
already published for $C^0$ finite element spaces for the
multi-frontal solver \cite{Calo2011,Paszynski2008}.

The restriction of this work to structured grid meshes is due to the
target application being B-spline-based isogeometric finite
elements. Isogeometric analysis is a relatively new method that has
been the subject of much work in recent years. While at its core it is
spline-based isoparametric finite element analysis, the spaces used
possess unique refinement strategies that form spaces that are
supersets of conventional finite elements. The mesh in B-spline-based
isogeometric analysis is always a structured grid of uniform
polynomial order. This simplification allows for specialized
complexity and memory usage estimates to be developed for these spaces
when the direct solver is the multi-frontal solver.

The higher continuous basis is known to approximate smooth functions
(e.g. solutions to PDEs) with orders of magnitude fewer degrees of
freedom than their $C^0$ counterparts. While this is a promising
result, the linear systems resulting from inner products of the higher
continuous basis functions are also orders of magnitude more expensive
to solve. This is something we have addressed in a previous paper
\cite{Collier2011b}. In this paper we explain the derivation of the
estimates in more detail and extend them to all spatial dimensions.

First we explain the concepts behind the multifrontal direct solver
algorithm and highlight why the algorithm is well-suited to handle
linear systems resulting from traditional $C^0$ finite elements
spaces. We emphasize the perspective that the algorithm may be viewed
as early $LU$-factorization. Second we evaluate computational
complexity of a single level of the multifrontal algorithm. Last we
generalize the single level result to all levels and provide estimates
for $C^0$ and $C^{p-1}$ B-splines spaces.

%% file: alg.tex
The state of the art direct solver for sparse linear systems is the
multi-frontal solver proposed by \cite{Duff1983,Duff1984}. It is the
generalization of the frontal solver algorithm proposed by
\cite{Irons1970}. Note that while in \cite{Irons1970}, the original
idea was developed for finite elements, the generalization applies to
general sparse linear systems.

The key observation in the algorithm is that $LU$-factorization may be
started during assembly on portions of the linear system that are
fully assembled on a local level. The general algorithm first examines
the matrix sparsity pattern to locate fully assembled matrix blocks
which are loosely connected to the remaining part of the matrix. When
used with finite elements, these blocks may be automatically
determined using what is known about the support of the basis
functions.

In the context of finite elements, the multifrontal direct solver
algorithm works in the following manner. We create the elemental
matrices using standard finite element procedures. In the direct
solver terminology, these elemental matrices are known as {\em frontal
  matrices}. The element matrices are typically assembled into a
global matrix, where contributions from shared degrees of freedom with
other elements are combined. However, depending on the topology and
order of the finite elements, there are degrees of freedom which at
the element level are fully assembled. For efficiency, these fully
assembled degrees of freedom are eliminted in terms of the partially
assembled ones at the element level, that is, at the frontal
level. Since this procedure can be repeated concurrently at each
element, it is usually known as {\em multi-frontal}. Particularizing
this idea to finite elements, we will start with an example of a two
element finite element mesh in any spatial dimension.
\subsection{Single level example: two element mesh}\label{s:sc}
Consider the partitioning of the elemental matrices in
equation~(\ref{eqn:e1}). We reorder the elemental matrices by first
listing those degrees of freedom that are fully assembled on the
element level, $x_e$, followed by those that are shared with other
elements, $y_e$, where subscript $e$ refers to the element
number. Thus the element matrix can be blocked accordingly to
represent interactions between fully assembled degrees of freedom and
those shared with other elements. Note that $A_e$ represents the block
of interactions which are fully assembled at the element level, blocks
$B_e$ and $C_e$ represent the interactions of fully assembled and
shared degrees of freedom, and block $D_e$ represents interactions of
shared degrees of freedom. Particularizing this for the two element
mesh, we obtain
\begin{equation}
\begin{bmatrix}
A_1 & B_1\\
C_1 & D_1
\end{bmatrix}\cdot\begin{bmatrix}
x_1\\
y_1
\end{bmatrix} = \begin{bmatrix}
f_1\\
g_1
\end{bmatrix},\hspace{0.5in}\begin{bmatrix}
A_2 & B_2\\
C_2 & D_2
\end{bmatrix}\cdot\begin{bmatrix}
x_2\\
y_2
\end{bmatrix} = \begin{bmatrix}
f_2\\
g_2
\end{bmatrix}\label{eqn:e1}
\end{equation}
Because the block $A_e$ is fully assembled, we may begin the
$LU$-factorization early and at the element level. Thus for each
element, we can multiply the top row by $C_eA_e^{-1}$ and subtract
from the bottom row,
\begin{equation}
\begin{bmatrix}
A_e & B_e\\
0 & D_e-C_eA_e^{-1}B_e
\end{bmatrix}\cdot\begin{bmatrix}
x_e\\
y_e
\end{bmatrix} = \begin{bmatrix}
f_e\\
g_e-C_eA_e^{-1}B_e
\end{bmatrix}\label{eqn:sc}
\end{equation}
Then, for each element $e=1,2$, the block matrix $D_e-C_eA_e^{-1}B_e$
and the vector $g_e-C_eA_e^{-1}B_e$ is assembled. After the
contributions from both elements are assembled, we can solve for
$y$. Once $y$ is computed, we can resort to backward substitution at
the element level to compute $x_e$.
\subsection{Multilevel example: eight element mesh in three dimensions}
\begin{figure}[ht]
\centering
\includegraphics[width=0.8\textwidth]{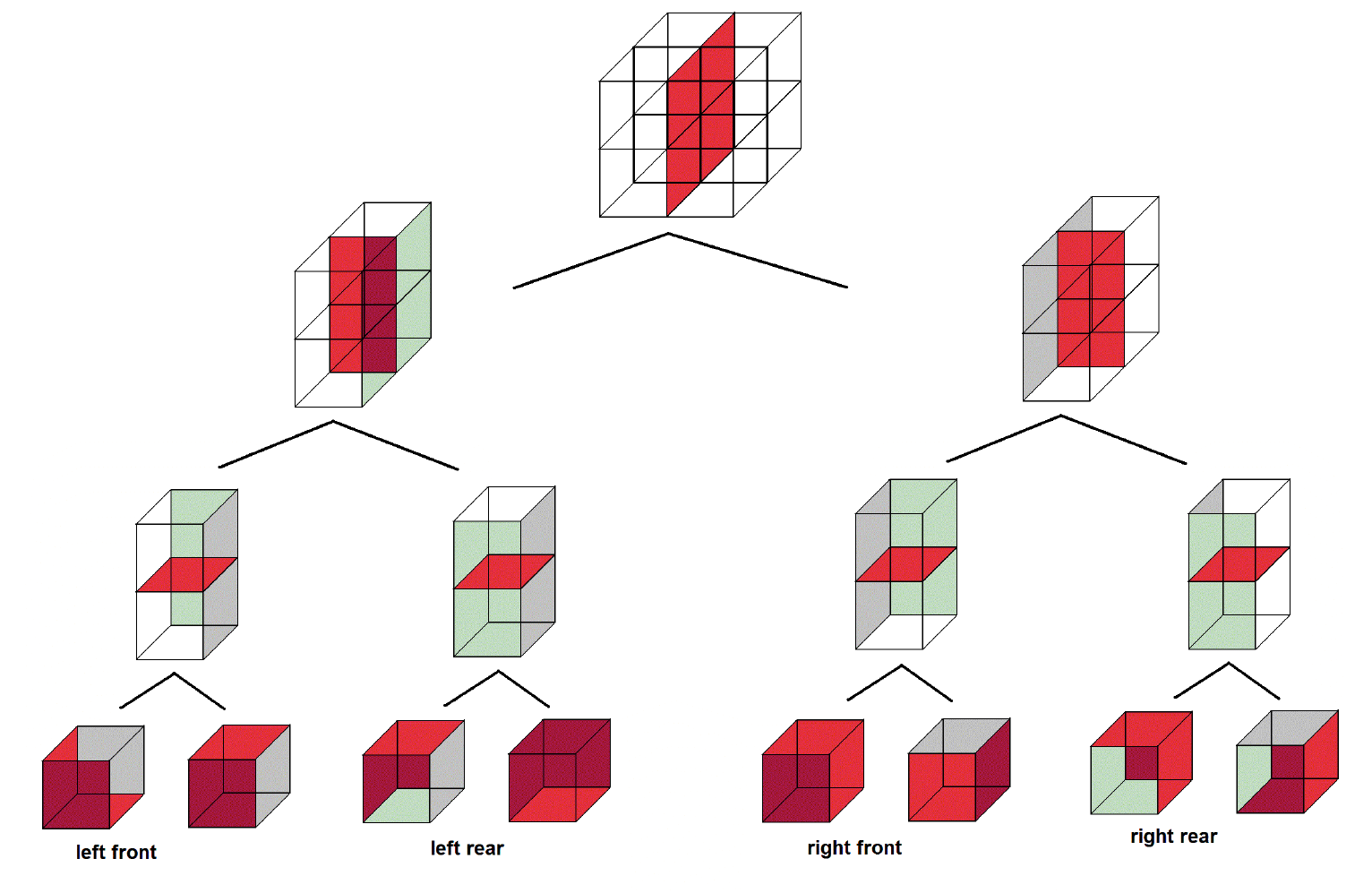}
\caption{The four levels of the elimination tree for a cube-shaped mesh with eight finite elements}\label{f:multifrontal}
\end{figure}

The procedure can be recursively generalized into multiple levels. For
example, the procedure for an eight element mesh is shown in
figure~\ref{f:multifrontal}. The elimination proceeds as follows:
\begin{enumerate}
\item Perform the local elimination of fully assembled degrees of
  freedom in each element as described in section~\ref{s:sc}. Note
  that the degrees of freedom eliminated, if any, are those which have
  support only on the element, the so-called bubble functions.
\item Pair the eight elements into any 4 clusters where the pairs of
  elements share a common face. To these frontal matrices, we apply
  the algorithm again, that is, we eliminate the fully assembled
  degrees of freedom in terms of the remaining degrees of freedom
  shared with other elements. At this level, the degrees of freedom
  eliminated are those with support on the shared face.
\item At the next level, we pair the four element clusters again into
  two which share a common interface. The recursive elimination
  procedure is applied at this level and repeated until we obtain a
  single cluster whose degrees of freedom are fully assembled (the top
  level in figure~\ref{f:multifrontal}).
\end{enumerate}
The connectivity graph describing the order of elimination and
clustering is called the {\em elimination tree}. At this point, with
the solution to the fully assembled system, we can move down the
elimintation tree, using backward substitution to recover the
remaining unknown degrees of freedom (those that were fully assembled
at each elimination level).

%% file: complexity.tex
Now we look at complexity comparing higher-continuous B-spline spaces
to their $C^0$ counterparts. For simplicity we will only consider
spaces which are $C^{p-1}$ continuous, that is, spaces which have
$p-1$ continuous derivatives across element interfaces, where $p$
refers to the polynomial order. Numerical tests indicate that these
spaces form limiting cases in the performance of the direct solver
algorithm.

The examples in the previous section come from $C^0$ finite element
spaces. In these examples, $(p-1)^d$ degrees of freedom may be
eliminated at the first level, where $d$ refers to the spatial
dimension. For $C^{p-1}$ spaces, the supports of the basis functions
spread into multiple elements and thus no degrees of freedom may be
eliminated at the element level. To be able to eliminate a degree of
freedom, we will cluster $p+1$ elements in each spatial dimension
together and use this on the first level for $C^{p-1}$ spaces. This is
the smallest grouping that can be obtained where a degree of freedom
is fully assembled.

In the remaining portion of this section we will develop the estimates
for complexity and memory usage. Computational complexity will be
estimated by counting floating point operations (FLOPS) and the memory
usage will we counted as the bytes needed to store the $LU$
factorization at each level of the elimination tree. The main building
block of the algorithm is the early $LU$-factorization, also known as
the Schur complement or static condensation. We will first develop the
cost of the Schur complement, and then proceed to apply this to the
full algorithm.
\subsection{Cost of the Schur complement}
In order to estimate the FLOPS and memory required to perform the
above partial $LU$ factorization, we will count operations and memory
used in forming the Schur complement, as shown in
equation~(\ref{eqn:sc}).  We denote the dimension of the square matrix
$A$ by $q$. We denote the number of columns in $B$ and rows in $C$ as
$r$, where $r$ is an assumed a constant. Then, we have:
\begin{equation}
\begin{tabular}{l}
  \mbox{FLOPS} = ${\cal O}( q^3 + q^2 r + q r^2) = {\cal O}( q^3 + q r^2)$\\[0.03in]
  \mbox{Memory}= ${\cal O}( q^2 + q r)$
\end{tabular}\label{eqn:schur}
\end{equation}
The FLOPS estimate is obtained by counting the operations needed to
find the $LU$ factors of $A$, $\mathcal{O}(q^3)$. To this we add the
FLOPS required to perform $r$ back-substitutions to form $A^{-1}B$,
$\mathcal{O}(rq^2)$. Finally, we add the cost of matrix multiplication
of $C$ to $A^{-1}B$, $\mathcal{O}(q r^2)$. The memory estimate is
obtained by adding the memory needed to store the $LU$ factors of the
matrix $A$, $\mathcal{O}(q^2)$, to that required to store $A^{-1}B$
and $CA^{-1}$, $\mathcal{O}(rq)$.

In the above memory estimate, we are only concerned with the space
required to store $L$ and $U$, since it is well-known that the cost of
storing original matrix is always smaller or equal than the memory
required to store factors $L$ and $U$. In particular, we have not
included the memory required to store the Schur complement, since this
is replaced in the next steps of $LU$ factorization by additional
Schur complement operations.
\subsection{Cost of the multi-frontal solver}
We divide our computational domain in $N_c$ clusters of elements. For
the $C^0$ case, each cluster is simply an element, while for
$C^{p-1}$, each cluster is a set of $p+1$ consecutive elements in each
spatial dimension. We assume for simplicity that the number of
clusters in our computational domain is $(2^d)^s$, where $s$ is a
positive integer which represents the number of levels of the
multi-frontal algorithm. Notice that even if this assumption is not
verified, the final result still holds true provided that the number
of degrees of freedom is sufficiently large.

The multi-frontal direct solver algorithm is summarized in algorithm
~\ref{a:multifrontal}. The FLOPS and memory required by
algorithm~\ref{a:multifrontal} can be expressed as
\begin{equation}
     \sum_{i=0}^{s-1} N_c(i)S(i)\label{eqn:cost}
\end{equation}
where $S(i)$ is the cost (either FLOPS or memory) of performing each
Schur complement at the $i^{th}$ level. Using the notation of the
previous subsection on the Schur complement, we define $q=q(i)$ as the
number of interior unknowns of each cluster at the $i^{th}$ step, and
$r=r(i)$ as the number of interacting unknowns at the $i^{th}$
step. We construct estimates for these numbers and summarize them in
table~\ref{tab:q_r_estimates}.
\begin{algorithm}
\caption{Multi-Frontal Algorithm\label{a:multifrontal}}
\begin{algorithmic}[1]
  \For{$i=0 \mbox{ to } s-1$}
  \State $N_c=N_c(i)=(2^d)^{s-i}$
  \If{$i=0$}
  \State Define $N_c(0)$ clusters
  \Else
  \State Join the old $N_c(i-1)$ clusters
  \State Eliminate interior degrees of freedom
  \State Define $N_c(i)$ new clusters
  \EndIf
  \EndFor
\end{algorithmic}
\end{algorithm}
\begin{table}[htp]
\centering
\caption{\label{tab:q_r_estimates} Number of interior ($q$) and interacting ($r$)
  unknowns at each level $i$ of the multi-frontal solver.}
\begin{tabular}{lcccccc}
\hline 
& & $q(0)$ & $r(0)$ & & $q(i),\ i \neq 0$  &  $r(i),\ i \neq 0$ \\
\hline
$C^0$  & & ${\cal O}(p^d)$ &  ${\cal O}( p^{d-1})$ & &
${\cal O}( 2^{(d-1)i} p^{d-1})$ &  ${\cal O}( 2^{(d-1)i} p^{d-1})$ \\
$C^{p-1}$ & & ${\cal O}( 1)$  & ${\cal O}(p^d)$ & &
${\cal O}(2^{(d-1)i} p^d)$ & ${\cal O}( 2^{(d-1)i}  p^d)$ \\
\hline
\end{tabular}
\end{table}

Let $N$ be the total number of unknowns in the original system. We use
the results from table~\ref{tab:q_r_estimates} with the FLOPS and
memory estimates in equation~(\ref{eqn:schur}) to develop
table~\ref{tab:S_i_estimates}. This table describes the cost in FLOPS
and memory of each level of the multi-frontal algorithm.
\begin{table}[htp]
\centering
  \caption{\label{tab:S_i_estimates}FLOPS and memory estimates at each
    level $i$ of the multi-frontal solver.}
  \begin{tabular}{lcccccc}
    \hline 
    & & FLOPS & Memory & & FLOPS & Memory\\
    & &      S(0)       &  S(0)    & &
    S(i) , $i \neq 0$  &  S(i) , $i \neq 0$      \\
    \hline
    $C^0$  & & ${\cal O}(p^{9})$ & ${\cal O}(p^{6})$ & & ${\cal O}(
    2^{6i} p^{6})$ & ${\cal O}( 2^{4i} p^{6})$   \\
    $C^{p-1}$ & & ${\cal O}( p^{6})$  & ${\cal O}(p^3)$ & & ${\cal O}(
    2^{6i} p^{9})$ & ${\cal O}( 2^{4i} p^{6})$    \\
    \hline
  \end{tabular}
\end{table}

Finally we use equation~(\ref{eqn:cost}) and table~\ref{tab:S_i_estimates}
to specialize estimates for $C^0$ and $C^{p-1}$ B-splines in one to
three spatial dimensions.
\paragraph{Estimates for 1D $C^0$ B-splines}
\begin{equation*}
\begin{tabular}{ll}
FLOPS = & $\ds 2^{s} p^3 + \sum_{i=1}^{s-1} 2^{s-i} = {\cal O}(2^s p^3) = {\cal
  O}(N_p p^3) = {\cal O}(N p^2)$\\
Memory = & $\ds 2^{s} p^2 + \sum_{i=1}^{s-1} 2^{s-i} = {\cal O}(2^s p^2) = {\cal
  O}(N_p p^2) = {\cal O}(N p)$
\end{tabular}
\end{equation*}
\paragraph{Estimates for 1D $C^{p-1}$ B-splines}
\begin{equation*}
\begin{tabular}{ll}
FLOPS = & $\ds 2^{s} p^2 + \sum_{i=1}^{s-1} 2^{s-i} p^3 = {\cal O}(2^s p^3) = {\cal
  O}(N_p p^3) = {\cal O}(N p^2)$,\\
Memory = & $\ds 2^{s} p + \sum_{i=1}^{s-1} 2^{s-i} p^2 = {\cal O}(2^s p^2) = {\cal
  O}(N_p p^2) = {\cal O}(N p)$.
\end{tabular}
\end{equation*}
\paragraph{Estimates for 2D $C^0$ B-splines}
\begin{equation*}
\begin{tabular}{ll}
FLOPS =& $\ds 2^{2s} p^6 + \sum_{i=1}^{s-1} 2^{2(s-i)} 2^{3i} p^3 = {\cal
  O}(2^{2s} p^6 + 2^{3s} p^3) =$\\
& $\ds {\cal O}(N_p^2 p^6 + N_p^3 p^3) =
{\cal O}(N p^4 + N^{1.5})$\\
Memory =& $\ds 2^{2s} p^4 + \sum_{i=1}^{s-1} 2^{2(s-i)} 2^{2i} p^2 = {\cal
  O}(2^{2s} p^4 + s^2 2^{2s} p^2) =$\\
& $\ds {\cal O}(N_p^2 p^4 + N_p^2 p^2 \log (N_p^2/p^2)) = {\cal O}(N p^2 + N \log (N/p^2))$
\end{tabular}
\end{equation*}
\paragraph{Estimates for 2D $C^{p-1}$ B-splines}
\begin{equation*}
\begin{tabular}{ll}
FLOPS =& $\ds 2^{2s} p^4 + \sum_{i=1}^{s-1} 2^{2(s-i)} 2^{3i} p^6 = {\cal
  O}(2^{2s} p^4 + 2^{3s} p^6) =$\\
& $\ds {\cal O}(N_p^3 p^6) = {\cal O}(N^{1.5} p^3)$\\
Memory =& $\ds 2^{2s} p^2 + \sum_{i=1}^{s-1} 2^{2(s-i)} 2^{2i} p^4 = {\cal
  O}(2^{2s} p^2 + s^2 2^{2s} p^4) =$\\
& $\ds {\cal O}(N_p^2 p^4 \log (N_p^2/p^2)) = {\cal O}(p^2 N \log (N/p^2))$
\end{tabular}
\end{equation*}
\paragraph{Estimates for 3D $C^{0}$ B-splines}
\begin{equation*}
\begin{tabular}{ll}
FLOPS =& $\ds 2^{3s} p^9 + \sum_{i=1}^{s-1} 2^{3(s-i)} 2^{6i} p^6 = {\cal
  O}(2^{3s} p^9 + 2^{6s} p^6) =$\\
& $\ds {\cal O}(N_p^3 p^9 + N_p^6 p^6) =
{\cal O}(N p^6 + N^2)$\\
Memory =&$\ds 2^{3s} p^6 + \sum_{i=1}^{s-1} 2^{3(s-i)} 2^{4i} p^4 = {\cal
  O}(2^{3s} p^6 + 2^{4s} p^4) =$\\
& $\ds {\cal O}(N_p^3 p^6 + N_p^4 p^4) =
{\cal O}(N p^3 + N^{4/3})$
\end{tabular}
\end{equation*}
\paragraph{Estimates for 3D $C^{p-1}$ B-splines}
\begin{equation*}
\begin{tabular}{ll}
FLOPS =& $\ds 2^{3s} p^6 + \sum_{i=1}^{s-1} 2^{3(s-i)} 2^{6i} p^9 = {\cal
  O}(2^{3s} p^6 + 2^{6s} p^9) =$\\
& $\ds {\cal O}(N_p^3 p^6 + N_p^6 p^9) =
{\cal O}(N^2 p^3)$\\
Memory =&$\ds 2^{3s} p^4 + \sum_{i=1}^{s-1} 2^{3(s-i)} 2^{4i} p^6 = {\cal
  O}(2^{3s} p^4 + 2^{4s} p^6) =$\\
& $\ds {\cal O}(N_p^4 p^6) =
{\cal O}(p^2 N^{4/3})$ 
\end{tabular}
\end{equation*}

%% file: results.tex
To test the validity of these estimates, we compute solutions to the
Laplace equation in three spatial dimensions on the unit cube
\begin{equation}
\begin{cases}
  -\nabla\cdot(\nabla u)=0&\ \ \ \text{on }\Omega\\
  u=0&\ \ \ \text{on }\Gamma_{D0}\\
  u=1&\ \ \ \text{on }\Gamma_{D1}\\
  (\nabla u)\cdot\mathbf{n}=0&\ \ \ \text{on }\Gamma_N
\end{cases}
\end{equation}
where $\Omega=[0,1]^3$, $\Gamma_{D0}=(:,:,0)$, $\Gamma_{D1}=(:,:,1)$,
and $\Gamma_N=(0,:,:) \cup (1,:,:) \cup (:,0,:) \cup (:,1,0)$.

All computational experiments have been performed on a workstation
with two quad-core Xeon X5550 processors and 24 Gb of memory running
Fedora 11. The model problem was implemented using PETSc
\cite{petsc1,petsc2} data structures. We used the MUMPS
\cite{Amestoy2001,Amestoy2006} implementation of the multi-frontal
algorithm, with METIS \cite{METIS} ordering (nested dissection). We
interfaced to MUMPS through PETSc, with the option to solve an
asymmetric system. Note that only one core was used in these numerical
experiments.

For all numerical results, we relate the FLOPS estimates to the
computational time measured for the solution of the linear system. The
memory estimates we relate to the number of nonzero entries in the
$LU$ factors. We report the memory required to store an integer and a
double precision number for each nonzero entry. Note that this is a
conservative quantification of memory usage. In general, solvers will
require the use of additional memory. However, we chose to report the
memory required to store the $LU$ factors as it is most closely
related to the estimates derived in the previous section.

We tested the estimates by solving the model problem for a linear
system containing 100,000 degrees of freedom for both $C^0$ and
$C^{p-1}$ spaces. This was accomplished by using different numbers of
elements. Note that while this will affect assembly time, this is not
included here in the time reported.

\begin{figure}[ht]
  \centering
  \subfloat[$C^0$ B-spline spaces]{\includegraphics[width=0.49\textwidth]{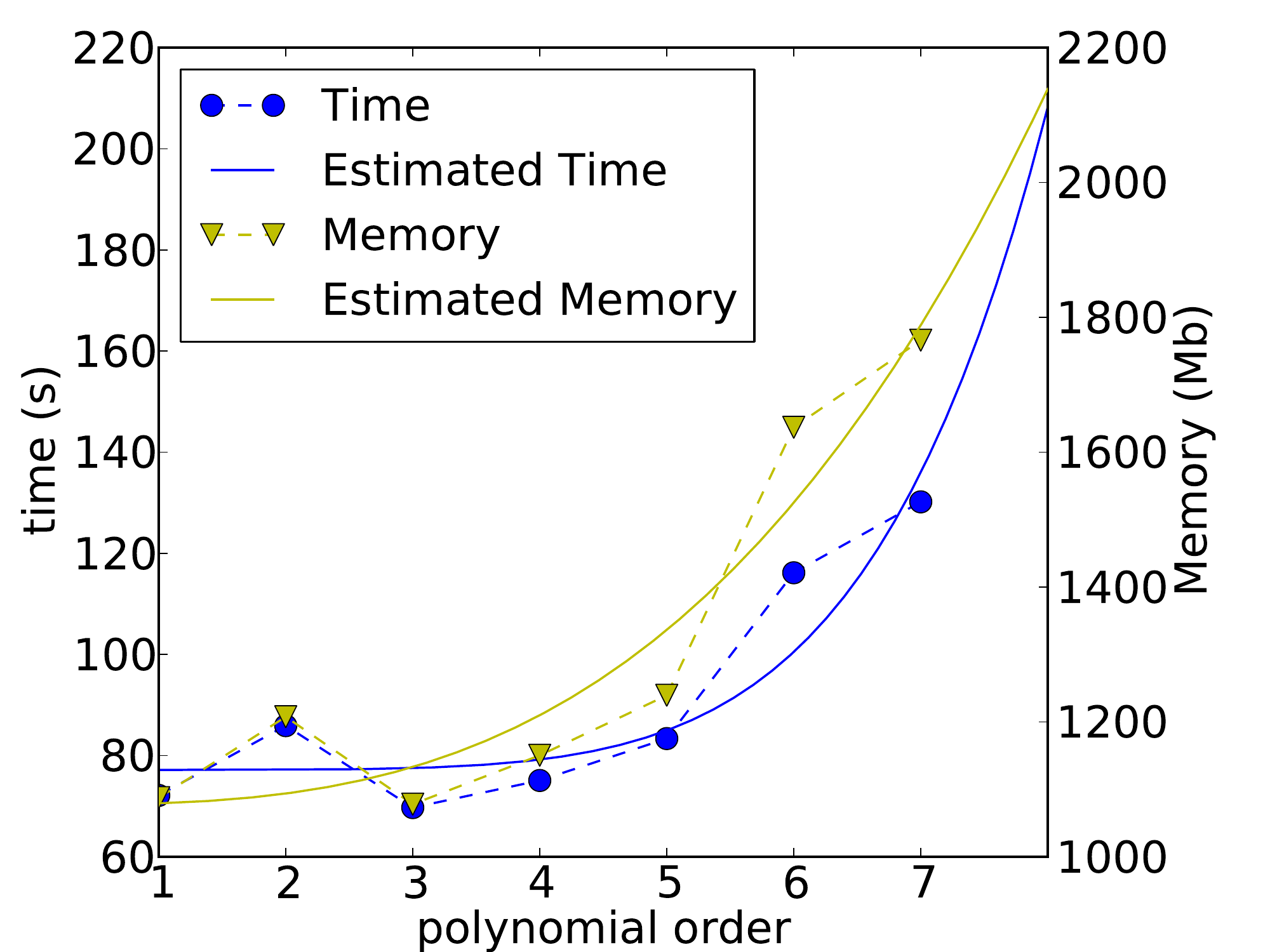}}
  \subfloat[$C^{p-1}$ B-spline spaces]{\includegraphics[width=0.49\textwidth]{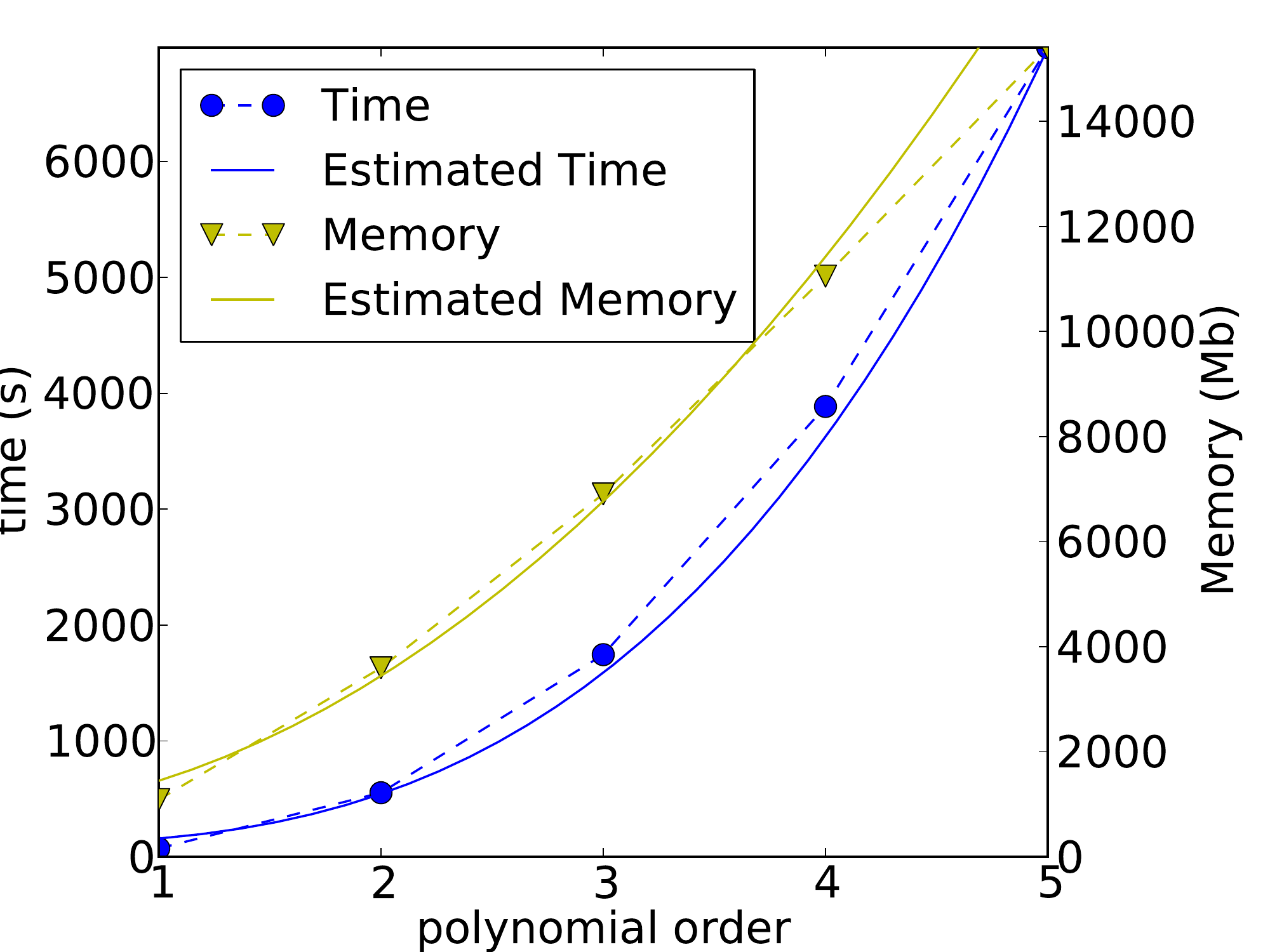}}
\caption{Time and memory for 100,000 degrees of freedom systems along with estimates}\label{f:valid}
\end{figure}

The numeric results of this test are shown in figure \ref{f:valid}. We
show computational time and memory usage for the $C^0$ spaces. The
solid line shown in a statistical best fit of the data to the
estimate. The excellent agreement between the estimates and the real
data supports the accuracy of the estimates.

%% file: conc.tex
In this paper, we derive estimates for the computational complexity
and memory usage of the multifrontal direct solver algorithm applied
to finite elements, where the mesh is a structured grid. This
restriction is with a view to forming estimates for higher continuous
spaces, $C^{p-1}$ B-spline spaces. We present numerical results to
support the validity of these estimates.

We observe that the estimates for 1D, 2D, and 3D are essentially
different. This is because of the structure of what can be eliminated
at each level. Although in 1D the number of FLOPS is linear with
respect the number of unknowns $N$, this case is not interesting. We
note that for $C^0$ spaces in two and three spatial dimensions, the
number of FLOPS is independent of $p$ if the number of unknowns $N$ is
large enough.  Therefore, under the assumption that $N$ is large
enough, an adaptive algorithm should select refinements based {\em
  exclusively} on the maximum decrease of error per added unknown,
independently of wether the refinement takes place on $h$ or $p$,
where $h$ here refers to the support of the basis function. This
efficiency is exploited in the adaptive algorithms used in $hp$-finite
elements \cite{hpfem}.

For 2D and 3D, the number of FLOPS of the $C^{p-1}$ spaces is $p^3$
times more expensive than the $C^0$ spaces, provided $N$ is large
enough. Thus, an adaptive algorithm for isogeometric analysis {\em
  should not} be based exclusively on the maximum decrease of the
error per added unknown. It would need to incorporate a special
treatment of the cost of each added unknown depending upon the type of
refinement.

%% file: Estimates.bbl
\begin{thebibliography}{10}
\expandafter\ifx\csname url\endcsname\relax
  \def\url#1{\texttt{#1}}\fi
\expandafter\ifx\csname urlprefix\endcsname\relax\def\urlprefix{URL }\fi
\expandafter\ifx\csname href\endcsname\relax
  \def\href#1#2{#2} \def\path#1{#1}\fi

\bibitem{Farin2002}
G.~Farin, Curves and Surfaces for CAGD: A Practical Guide, 5th Edition, Morgan
  Kaufmann, 2002.

\bibitem{Piegl1995}
L.~Piegl, W.~Tiller, The NURBS Book, Monographs in Visual Communication,
  Springer, New York, 1995.

\bibitem{Hughes2005}
T.~J.~R. Hughes, J.~Cottrell, Y.~Bazilevs, Isogeometric analysis: {CAD}, finite
  elements, {NURBS}, exact geometry and mesh refinement, Computer Methods in
  Applied Mechanics and Engineering 194 (2005) 4135--4195.

\bibitem{Cottrell2006}
J.~Cottrell, A.~Reali, Y.~Bazilevs, T.~J.~R. Hughes, Isogeometric analysis of
  structural vibrations., Computer Methods in Applied Mechanics and Engineering
  195~(41-43) (2006) 5257.

\bibitem{Calo2011}
V.~M. Calo, N.~O. Collier, D.~Pardo, M.~R. Paszynski,
  \href{http://www.sciencedirect.com/science/article/pii/S1877050911002596}{Computational
  complexity and memory usage for multi-frontal direct solvers used in p finite
  element analysis}, Procedia Computer Science 4 (2011) 1854 -- 1861,
  proceedings of the International Conference on Computational Science, ICCS
  2011.
\newblock \href {http://dx.doi.org/DOI: 10.1016/j.procs.2011.04.201}
  {\path{doi:DOI: 10.1016/j.procs.2011.04.201}}.
\newline\urlprefix\url{http://www.sciencedirect.com/science/article/pii/S1877050911002596}

\bibitem{Paszynski2008}
M.~Paszynski, Performance of multi level parallel direct solver for hp finite
  element method, in: R.~Wyrzykowski, J.~Dongarra, K.~Karczewski, J.~Wasniewski
  (Eds.), Parallel Processing and Applied Mathematics, Vol. 4967 of Lecture
  Notes in Computer Science, Springer Berlin / Heidelberg, 2008, pp.
  1303--1312.

\bibitem{Collier2011b}
N.~Collier, D.~Pardo, L.~Dalcin, M.~Paszynski, V.~M. Calo, The cost of
  continuity: a study of the performance of isogeometric finite elements using
  direct solvers, submitted to Computer Methods in Applied Mechanics and
  Engineering.

\bibitem{Duff1983}
I.~S. Duff, J.~K. Reid, The multifrontal solution of indefinite sparse
  symmetric linear, ACM Trans. Math. Softw. 9 (1983) 302--325.
\newblock \href {http://dx.doi.org/http://doi.acm.org/10.1145/356044.356047}
  {\path{doi:http://doi.acm.org/10.1145/356044.356047}}.

\bibitem{Duff1984}
I.~S. Duff, J.~K. Reid, The multifrontal solution of unsymmetric sets of linear
  equations, SIAM Journal on Scientific and Statistical Computing 5~(3) (1984)
  633--641.

\bibitem{Irons1970}
B.~M. Irons, A frontal solution program for finite element analysis,
  International Journal for Numerical Methods in Engineering 2~(1) (1970)
  5--32.
\newblock \href {http://dx.doi.org/10.1002/nme.1620020104}
  {\path{doi:10.1002/nme.1620020104}}.

\bibitem{petsc1}
S.~Balay, K.~Buschelman, W.~D. Gropp, D.~Kaushik, M.~G. Knepley, L.~C. McInnes,
  B.~F. Smith, H.~Zhang, {PETSc} {W}eb page, http://www.mcs.anl.gov/petsc
  (2010).

\bibitem{petsc2}
S.~Balay, K.~Buschelman, V.~Eijkhout, W.~D. Gropp, D.~Kaushik, M.~G. Knepley,
  L.~C. McInnes, B.~F. Smith, H.~Zhang, {PETS}c users manual, Tech. Rep.
  ANL-95/11 - Revision 3.0.0, Argonne National Laboratory (2008).

\bibitem{Amestoy2001}
P.~R. Amestoy, I.~S. Duff, J.~Koster, J.-Y. L'Excellent, A fully asynchronous
  multifrontal solver using distributed dynamic scheduling, {SIAM} Journal of
  Matrix Analysis and Applications 23~(1) (2001) 15--41.

\bibitem{Amestoy2006}
P.~R. Amestoy, A.~Guermouche, J.-Y. L'Excellent, S.~Pralet, Hybrid scheduling
  for the parallel solution of linear systems, Parallel Computing 32~(2) (2006)
  136--156.

\bibitem{METIS}
G.~Karypis, V.~Kumar, Parallel multilevel k-way partitioning scheme for
  irregular graphs, in: Proceedings of the 1996 ACM/IEEE Conference on
  Supercomputing, 1996, p.~35.

\bibitem{hpfem}
L.~Demkowicz, Computing with $hp$-adaptive finite elements, vol. 1: One and two
  dimensional elliptic and {M}axwell problems, Chapman \& Hall/CRC, 2007.

\end{thebibliography}


18 gid=1840429327
18 uid=2022889333
20 ctime=1333870967
20 atime=1333871770
24 SCHILY.dev=234881026
22 SCHILY.ino=8034429
18 SCHILY.nlink=1
